\documentclass[reqno,twoside,12pt]{amsart}

\hoffset - 1.5cm

\usepackage{graphicx}
\usepackage{amstext}
\usepackage{pstricks,pst-node,pst-text,pst-3d}
\usepackage{amsmath}
\usepackage{amsthm}
\usepackage{amssymb}

\newtheorem{Theorem}{Theorem}[section]
\newtheorem{Fact}{Fact}[section]
\newtheorem{Lemma}{Lemma}[section]
\newtheorem{Proposition}{Proposition}[section]
\newtheorem{Corollary}{Corollary}[section]
\theoremstyle{definition}
\newtheorem{Definition}{Definition}[section]
\newtheorem{Example}{Example}[section]
\newtheorem{Remark}{Remark}[section]
\newtheorem{Hipothesis}{Hipothesis}[section]

\newcommand{\ba}{\begin{array}}
\newcommand{\bc}{\begin{center}}
\newcommand{\bd}{\begin{description}}
\newcommand{\bdm}{\begin{displaymath}}
\newcommand{\be}{\begin{enumerate}}
\newcommand{\beq}{\begin{equation}}
\newcommand{\bdf}{\begin{Definition}}
\newcommand{\bex}{\begin{Example}}
\newcommand{\bft}{\begin{Fact}}
\newcommand{\bl}{\begin{Lemma}}
\newcommand{\bp}{\begin{Proposition}}
\newcommand{\br}{\begin{Remark}}
\newcommand{\bt}{\begin{Theorem}}
\newcommand{\bco}{\begin{Corollary}}
\newcommand{\bh}{\begin{Hipothesis}}
\newcommand{\ea}{\end{array}}
\newcommand{\ec}{\end{center}}
\newcommand{\ed}{\end{description}}
\newcommand{\edm}{\end{displaymath}}
\newcommand{\ee}{\end{enumerate}}
\newcommand{\eeq}{\end{equation}}
\newcommand{\edf}{\end{Definition}}
\newcommand{\eex}{\end{Example}}
\newcommand{\eft}{\end{Fact}}
\newcommand{\el}{\end{Lemma}}
\newcommand{\ep}{\end{Proposition}}
\newcommand{\er}{\end{Remark}}
\newcommand{\et}{\end{Theorem}}
\newcommand{\eco}{\end{Corollary}}
\newcommand{\eh}{\end{Hipothesis}}

\newcommand{\R}{\ensuremath{\mathbb{R}}}

\newcommand{\CC}{\mathcal{C}}
\newcommand{\CF}{\ensuremath{\mathcal{F}}}

\newcommand{\cl}{\mbox{\normalfont{Cl}}}

\def\e{\varepsilon}
\begin{document}

\title[Periodic orbits for a class of galactic potentials]{Periodic orbits for a class of galactic potentials}

\author{Felipe Alfaro}
\address{Academia de Matem\'aticas, UACM, Plantel Casa Libertad, M\'exico, D.F., 09620 M\'exico}
\email{eltiopi@yahoo.com}

\author{Jaume Llibre}
\address{Departament de Matem\`aticues\\
Universitat Aut\`onoma de Barcelona\\
Bellaterra, Barcelona 08193, Catalonia, Spain}
\email{jllibre@mat.uab.cat}

\author{Ernesto P\'{e}rez-Chavela}
\address{Departamento de Matem\'{a}ticas\\
UAM-Iztapalapa \\ Av. Rafael Atlixco 186, M\'{e}xico, DF 09340\\ Mexico}
\email{epc@xanum.uam.mx}
 
\keywords{Periodic orbit, averaging theory, Hamiltonian system.}

\subjclass[2010]{Primary: 34C29; Secondary: 37C27.}

\thanks{All authors has been partially supported Conacyt M\'exico, grant 128790. The second author is partially
supported by MCYT/FEDER grant number MTM2005-06098-C02-01, by a CICYT grant number 2005SGR 00550 and
by ICREA Academia.}

\maketitle

\date{\today}

\begin{abstract}
In this work, applying general results from averaging theory, we find  periodic orbits for a class 
of Hamiltonian systems $H$ whose potential models the motion of elliptic galaxies.
\end{abstract}


\section{introduction and statement of the main results}
Galactic dynamics is a branch of Astrophysics whose development started  only around sixty years ago, when it was possible to have a view of the physical world beyond the integrable and near integrable systems \cite{Co}. Even the importance 
of the analysis of galactic potentials, the global dynamics of galaxies is not a simple question and represents a big
challenge for the researches in the field \cite{BT}. Most of the work in the analysis of galaxies is numerical, in this paper we 
present an analytical technique, the averaging theory, which allows to find periodic orbits of a differential 
system.

In the last years, great quantity of the research on galactic dynamics has been focused on models of elliptical galaxies. In most of these models the terms in the potential are of even order, so we have adopted this fact in the Hamiltonian system that we are analyzing. Another important point that appears in these kind of potentials is that the existence of periodic orbits 
is a useful tool for constructing new and more complicated self consistent models. One way to identify periodic orbits is to localize the central fixed points on the surfaces of constant energy.
In \cite{PZ}, the authors study the localization of periodic orbits and their linear stability for a particular two-component galactic potential. In fact, in our days the study of individual orbits in some galactic potentials is a new branch of galactic dynamics (see for instance the articles \cite{Ca, ES, G}) . 

The calculation of particular orbits in some analytical potentials modeling elliptical galaxies, indicates that relatively small symmetry breaking  corrections can increase dramatically the number of stochastic orbits, showing
the importance of the study of perturbations of simple models \cite{HKM}. The class of potentials studied in this paper have not chosen with the aim of modeling some particular galaxies, our objective is to study systems which are generic in their basic properties.

\medskip

In \cite{Puc}, the authors study the galactic potential 
$$H = \frac{1}{2}(P_X^2 + P_Y^2) + V(X^2, Y^2).$$
These kind of potentials are important in the modeling of elliptic galaxies, as for instance we can mention the potentials 
$V_L = \log{(1+X^2+Y^2/q)}$ and $V_C= \sqrt{1+X^2+Y^2/q} - 1$, where the parameter $q$ gives the eccentricity of 
the elliptic galaxy.  In this paper we deal with  the Hamiltonian
\beq\label{ham-nr}
H = \frac{1}{2}(P_X^2 + X^2) + 
  \frac{1}{2q}  (P_Y^2 + Y^2) + (a X^4  + b X^2 Y^2 + 
      c Y ^4),
\eeq     
and its respective Hamilton's equation
\begin{eqnarray}\label{eq:ham1}
\dot {X} &=& P_X, \nonumber \\ \dot {Y} &=& \frac{P_Y}{q}, \\
\dot {P_X} &=& -X - (4aX^3 + 2bXY^2 ), \nonumber \\ 
\dot {P_Y} &=& -\frac{Y}{q} - ( 2bX^2Y  + 4cY^3), \nonumber
\end{eqnarray}
the matrix of the linear part of this system at the origin of coordinates is 
$$M = \left(\begin{array}{cccc}
    0  &  0 & 1 & 0  \\
   0  &  0 & 0 & \frac{1}{q}  \\
    -1  &  0 & 0 & 0  \\
     0  &  - \frac{1}{q} & 0 & 0  
    \end{array}\right),$$
 with eigenvalues $\pm i, \pm i/q$ where $i=\sqrt{-1}$.
 
 In order to obtain periodic orbits for these kind of potentials we will apply averaging theory, in this way 
 we re-parametrize the coordinates by the factor $\sqrt{\varepsilon}$ for $\varepsilon$ positive small enough (a similar change of coordinates has been used in \cite{LR}), that is
 we do the change 
\beq\label{rep}
(X, Y, P_X, P_Y) \to \sqrt{\varepsilon}(x, y, p_x, p_y).
\eeq

After straightforward computations we get the Hamilton's equations
\begin{eqnarray}\label{eq:ham2}
\dot {x} &=& p_x, \nonumber \\ \dot {y} &=& \frac{p_y}{q}, \\
\dot {p_x} &=& -x - \varepsilon(4ax^3 + 2bxy^2 ), \nonumber \\ \dot {p_y} &=& -\frac{y}{q} -
\varepsilon( 2bx^2y  + 4cy^3), \nonumber
\end{eqnarray}
which have the same linear part at the origin than the previous one, the structure of the new Hamiltonian is
identity with (\ref{ham-nr}) in the new variables. Our goal is to study which periodic orbits for $\varepsilon = 0$ (the unperturbed system) persists  for $\varepsilon$ positive and small enough (the perturbed system).

By the form of the matrix $M$  we observe the necessity to split the analysis for the periodic orbits in two cases
 \begin{itemize}
 \item $q$ is an irrational number. Here the linear part of system (\ref{eq:ham2}) has two planes foliated by periodic orbits. In the first one the orbits have period $2\pi$, each periodic orbit on this plane is of the form $$PO_1=(x_0\cos{t} + p_{x_0} \sin{t}, 0, p_{x_0} \cos{t} - x_0\sin{t}, 0).$$
 In the second one, the orbits have period $2\pi q$, each periodic orbit on this plane is of the form 
 $$PO_2=(0, y_0\cos{(t/q)} + p_{y_0}\sin{(t/q)}, 0, p_{y_0}\cos{(t/q)} - y_0\sin{(t/q))}.$$
 \item $q$ is a rational number. Here the linear part of system (\ref{eq:ham2}) has a $4$--dimensional space filled of 
 periodic orbits of period $2\pi r$ if $q=r/s$ with $(r,s)=1$, where each periodic orbit is of the form
 \begin{eqnarray} 
 PO_3 &=& {\rm (}x_0\cos{t} + p_{x_0} \sin{t},  y_0\cos{(st/r)} + p_{y_0}\sin{(st/r)}, \nonumber \\
 & & p_{x_0} \cos{t} - x_0\sin{t},
 p_{y_0}\cos{((st/r)} - y_0\sin{(st/r)} {\rm )}.\nonumber
 \end{eqnarray}
\end{itemize}

When $q$ is an irrational number our main result is:

\bt\label{main1}
For $q$ an irrational number, we have that in every energy level  
$H=h>0$ the Hamiltonian system \eqref{eq:ham1}  has 
\begin{itemize}
\item  [{\bf (a)}] at least one periodic solution $(X(t),Y(t),P_X(t),P_Y(t)$ such that when $\varepsilon \to 0,$
we have that $(X(0),Y(0),P_X(0),P_Y(0))$ tends to  $(0,0,0,0);$
\item [{\bf (b)}] at least one periodic solution $(X(t),Y(t),P_X(t),P_Y(t))$ such that when $\varepsilon \to 0,$ 
we have that $(X(0),Y(0),P_X(0),P_Y(0))$ tends to  $(0,0,0,0).$ 
\end{itemize}
So, for $q$ irrational, we obtain that in every energy level $H = h > 0$ the perturbed Hamiltonian system has at least $2$ periodic orbits.
\et

\br We note that the periodic orbits found in the statements of Theorem \ref{main1} are in fact degenerate Hopf bifurcations periodic orbits, since they born from the equilibrium point localized at the origin of coordinates. 
Unfortunately we cannot obtain periodic solutions when $q$ is a rational number, see Remark
\ref{remarkf}. \er

The paper is organized as follows. In section $2$ we present  the theorem from averaging theory necessary to prove our main result. In section $3$ we give the proof of Theorem \ref{main1}.

\medskip

\section{Some results from averaging theory}
In order to have a self contained paper, in this section we present the basic results from the
averaging theory that are necessary for proving the main results of this  paper.

\medskip

\subsection{Results from averaging theory}
We consider the problem of the bifurcation of $T$--periodic
solutions from the differential system
\beq\label{eq:4}
{\bf x}^{\prime}(t)= F_0(t,{\bf x})+\e F_1(t,{\bf x}),
\eeq
where the functions
$F_0,F_1: \R \times \Omega \to \R^n$  are of class $\CC^2$ functions, $T$--periodic
in the first variable, and $\Omega $ is an open subset of $\R^n$. When $\e = 0$
we get the unperturbed system
\begin{equation}\label{eq:5}
{\bf x}^{\prime}(t)= F_0(t,{\bf x}).
\end{equation}
One of the main assumptions on the above system is that it 
has a submanifold of periodic solutions. A solution of system (\ref{eq:4}), for 
$\e$ sufficiently small  is given using the averaging theory. For a general
introduction to the averaging theory see the books of Sanders and
Verhulst \cite{SV}, and of Verhulst \cite{Ve}.

\medskip

Let ${\bf x}(t,{\bf z})$ be the solution of the unperturbed system
\eqref{eq:5} such that ${\bf x}(0,{\bf z})= {\bf z}$. We write
the linearization of the unperturbed system along the periodic
solution ${\bf x}(t,{\bf z})$ as
\begin{equation}\label{eq:6}
{\bf y}^{\prime}=D_{\bf x}{F_0}(t,{\bf x}(t,{\bf z})){\bf y}.
\end{equation}
In what follows we denote by $M_{\bf z}(t)$ some fundamental
matrix of the linear differential system \eqref{eq:6}, and by
$\xi:\R^k\times \R^{n-k}\to \R^k$ the projection of $\R^n$ onto
its first $k$ coordinates; i.e. $\xi(x_1,\ldots,x_n)=
(x_1,\ldots,x_k)$.

\medskip

\bt\label{t3}
Let $V\subset \R^k$ be open and bounded, and let $\beta_0:
\cl(V)\to \R^{n-k}$ be a $\CC^2$ function. We assume that
\begin{itemize}
\item[(i)] $\mathcal{Z}=\left\{ {\bf z}_{\alpha}=\left( \alpha,
\beta_0(\alpha)\right),~~\alpha\in \cl(V) \right\}\subset \Omega$
and that for each ${\bf z}_{\alpha}\in \mathcal{Z}$ the solution
${\bf x}(t,{\bf z}_{\alpha})$ of \eqref{eq:5} is $T$--periodic;

\item[(ii)] for each ${\bf z}_{\alpha}\in \mathcal{Z}$ there
is a fundamental matrix $M_{{\bf z}_{\alpha}}(t)$ of \eqref{eq:6}
such that the matrix $M_{{\bf z}_{\alpha}}^{-1}(0)- M_{{\bf
z}_{\alpha}}^{-1}(T)$ has in the upper right corner the $k\times
(n-k)$ zero matrix, and in the lower right corner a $(n-k)\times
(n-k)$ matrix $\Delta_{\alpha}$ with $\det(\Delta_{\alpha})\neq 0$.
\end{itemize}
We consider the function $\CF:\cl(V) \to \R^k$
\begin{equation}\label{eq:7}
\CF(\alpha)=\xi\left( \int _0^T M_{{\bf
z}_{\alpha}}^{-1}(t)F_1(t,{\bf x}(t,{\bf z}_{\alpha})) dt\right).
\end{equation}
If there exists $a\in V$ with $\CF(a)=0$ and $\displaystyle{\det
\left( \left({d\CF}/{d\alpha}\right)(a)\right)\neq 0}$, then there
is a $T$--periodic solution $\varphi (t,\e)$ of system
\eqref{eq:4} such that $\varphi(0,\e)\to {\bf z}_a$ as $\e\to 0$.
\et

Theorem~\ref{t3} goes back to Malkin \cite{Ma} and Roseau
\cite{Ro}, for a shorter proof see \cite{BFL}.

\section{Proof of the main Theorem}
In this section we give the proof of our main result Theorem \ref{main1}.

\subsection{Proof of Theorem \ref{main1}}
We know that the periodic orbits of a Hamiltonian system always appear in  cylinders
foliated by periodic orbits, each periodic orbit  corresponds to a different value of the energy $h$,
see for more details \cite{AMR}. In order to have
 isolated periodic orbits and be able to apply the averaging theory we 
fix the total energy $H=h$.  Computing $p_x$ in the energy level $H=h$ we get
\beq\label{moment-x}
p_x = \pm \sqrt{
   2 h - p_y^2/q - x^2 - y^2/q - \e(2 a x^4  +  
    2 b x^2y^2 +  2 c y^4 )}. 
\eeq
The fix value $h$ of the total energy  is determined by the initial periodic orbit, which in our case for the 
periodic orbit $PO_1$ it corresponds  to  $h= \frac{1}{2}(p_{x_0}^2 + x_0^2)$, choosing the sign $+$
for $p_x $, and expanding around $\e =0$ we obtain
\beq\label{px-expand}
p_x = \sqrt { p_{x_0}^2 - x^2 + x_0^2 - (y^2 + p_y^2)/q}
   - \e \frac{a x^4 + b x^2 y^2 + c y^4} {\sqrt { p_{x_0}^2 - x^2 + x_0^2 - (y^2 + p_y^2)/q}} .
\eeq

The equations of motion on the energy level $H = (p_{x_0}^2+x_0^2)/2$ are
\begin{eqnarray}\label{eq:ham-red}
\dot {x} &=& \sqrt { p_{x_0}^2 - x^2 + x_0^2 - (y^2 + p_y^2)/q}
   - \e \frac{a x^4 + b x^2 y^2 + c y^4} {\sqrt { p_{x_0}^2 - x^2 + x_0^2 - (y^2 + p_y^2)/q}}, \nonumber \\ 
\dot {y} &=& \frac{p_y}{q}, \\
 \dot {p_y} &=& -\frac{y}{q} - \e( 2bx^2y  + 4cy^3). \nonumber
\end{eqnarray}

In order to apply Theorem \ref{t3} to system (\ref{eq:ham-red}), let
\begin{equation}\label{ap-t3}
\begin{array}{l}
{\bf x}= (x,y,p_y),\\
F_0(t,{\bf x}) = \left(  \sqrt { p_{x_0}^2 - x^2 + x_0^2 - (y^2 + p_y^2)/q} , \quad p_y/q, \quad -y/q \right),\\
F_1(t,{\bf x}) =  \left( - \frac{a x^4 + b x^2 y^2 + c y^4}{\sqrt { p_{x_0}^2 - x^2 + x_0^2 - (y^2 + p_y^2)/q} }, \quad 0, \quad  -( 2bx^2y  + 4cy^3) \right).\\
\end{array}
\end{equation}
The set $\Omega = \{ (x,y,p_y) | q (p_{x_0}^2 - x^2 + x_0^2) - y^2 - p_y^2 \neq 0 \}$ is an open
subset of $\R^3$.  Clearly the above functions are of class 
$C^2(\Omega).$ The set $V$ of Theorem \ref{t3} is given by
$$V= \{ {\bf z} = (x_0,0,0) :  |x_0| < \rho \} \quad {\rm for\ some}\  \rho \ {\rm large\ enough.}$$

Let ${\bf x}(t,{\bf z})$ be the solution of the unperturbed system
\eqref{eq:5} such that ${\bf x}(0,{\bf z})= {\bf z}$. The variational equations of the unperturbed 
system along the periodic solution $PO_1$ are
\begin{equation}\label{variational}
{\bf y}^{\prime}=D_{\bf x}{F_0}(t,{\bf x}(t,{\bf z})){\bf y},
\end{equation}
where ${\bf y}$ is a $3 \times 3$ matrix.

The fundamental matrix $M(t)$ of the differential system (\ref{variational}) such that $M(0)$ is 
the identity matrix of $\R^3$ takes the simple form
\begin{equation}\label{var-mat}
M(t) = \left(\begin{array}{ccc}
    \cos{t} - x_0 \sin{(t)}/p_{x_0}  &  0  & 0  \\
   0  &  \cos{(t/q)} & \sin{(t/q)}   \\
      0 & - \sin{(t/q)} & \cos{(t/q)}   
    \end{array}\right),
\end{equation}  
    whose inverse is given by 
    
    $$M^{-1}(t) = \left(\begin{array}{ccc}
    p_{x_0}/( p_{x_0}\cos{t} - x_0 \sin{t}) &  0  & 0  \\
   0  &  \cos{(t/q)} & - \sin{(t/q)}   \\  0 & \sin{(t/q)} & \cos{(t/q)}   
    \end{array}\right).$$
  
 An easy computation shows that
 $$M^{-1}(0) - M^{-1}(2\pi) = \left(\begin{array}{ccc}
  0 &  0  & 0  \\
   0  &  2\sin^2{(\pi/q)} &  \sin{(2\pi/q)}   \\
      0 & - \sin{(2\pi/q)} & 2\sin^2{(\pi/q)}  
    \end{array}\right).$$
We observe that this matrix has a couple of zeros in the upper right corner of size $1 \times 2$;  the
determinant of the $ 2 \times 2$ matrix which appears in the lower right corner is $4\sin^2{(\pi/q)} \neq 0$    
because $q$ is an irrational number. Consequently all the assumptions of Theorem \ref{t3} are satisfied. Therefore 
we must compute the simple zeroes of the function $\CF$ defined in Theorem \ref{t3}. A straightforward computations shows that 
$$F_1(t,{\bf x}(t,{\bf z})) = \left( -\frac{a (x_0 \cos{t} + p_{x_0} \sin{t})^4}{p_{x_0}\cos{t} -x_0\sin{t}}, 0, 0 \right),$$
therefore we get
$$M^{-1}(t) F_1(t,{\bf x}(t,{\bf z}))  =  \left( -\frac{a p_{x_0}(x_0 \cos{t} + p_{x_0} \sin{t})^4}{( p_{x_0}\cos{t} 
-x_0\sin{t})^2}, 0, 0) \right).$$

Let $$f_1 (x_0) = \frac{1}{2 \pi} \int_0^{2\pi}  -\frac{a p_{x_0}(x_0 \cos{t} + p_{x_0}\sin{t})^4}{\left( p_{x_0}\cos{t}
 -x_0\sin{t}\right)^2} dt = 3a p_{x_0}( p_{x_0}^2 + x_0^2)/2 = 3ah \sqrt{2h - x_0^2}.$$
In the last equality we have gotten $p_{x_0}$  from the energy relation $h=( p_{x_0}^2 + x_0^2)/2$. So the solutions of
 $f_1=0$ are $x_0= \pm \sqrt{2h}$,  which are simple zeroes.  On the other hand we can verify that both zeroes generate the same periodic orbit. Then doing the rescaling (\ref{rep}) we obtain statement {\bf (a)} of Theorem \ref{main1}.

\medskip

For the proof of statement {\bf (b)}, as in the previous case we 
fix the value of the total energy  as  $H= (p_{x_0}^2 + x_0^2)/2q = h$ determined by 
the initial periodic orbit $PO_2$. Computing $p_y$ from the equation $H=h$ we obtain
\beq\label{moment-y}
p_y = \pm \sqrt{
   2 qh - qp_x^2 - qx^2 - y^2  - \e(2 aq x^4  +  
    2 b qx^2y^2 +  2 c qy^4 )}, 
\eeq
we choose the sign $+$ for $p_y$ and expand arounf $\e=0$ getting
\beq\label{py-expand}
p_y = \sqrt{
   2 qh - q(p_x^2 + x^2) - y^2}  - \e \frac{q( a x^4  +  
     b x^2y^2 +  2 y^4 )}{\sqrt{ 2 qh - q(p_x^2 + x^2) - y^2}}.
 \eeq

Now, we write the equations of motion on the energy level $H= (p_{x_0}^2 + x_0^2)/2q$
in the order $(y,x,p_y)$, they  are given by the system
\begin{eqnarray}\label{eq:ham-red-y}
\dot {y} &=& \frac{\sqrt{
   2 qh - q(p_x^2 + x^2) - y^2}}{q} - 
\e \frac{ a x^4  +  
     b x^2y^2 +  2 y^4 }{\sqrt{
   2 qh - q(p_x^2 + x^2) - y^2}}, \nonumber \\
\dot {x} &=& p_x, \\ 
\dot {p_y} &=& -\frac{y}{q} - \e( 2bx^2y  + 4cy^3). \nonumber
\end{eqnarray}

In order to apply Theorem \ref{t3} to system (\ref{eq:ham-red-y}) we are using the same notations
and definitions (with the obvious changes) than in the previous case.

Let ${\bf x}(t,{\bf z})$ be the solution of the unperturbed system
\eqref{eq:5} such that ${\bf x}(0,{\bf z})= {\bf z}$. The variational equations of the unperturbed 
system along the periodic solution $PO_2$ are
\begin{equation}\label{variational-y}
{\bf y}^{\prime}=D_{\bf x}{F_0}(t,{\bf x}(t,{\bf z})){\bf y},
\end{equation}
where ${\bf y}$ is a $3 \times 3$ matrix.

The fundamental matrix $M(t)$ of the differential system (\ref{variational-y}) such that $M(0)$ is 
the identity matrix of $\R^3$ takes the simple form
\begin{equation}\label{var-mat1}
M(t) = \left(\begin{array}{ccc}
    \cos{(t/q)} - y_0 \sin{(t/q)}/p_{y_0}  &  0  & 0  \\
   0  &  \cos{t} & \sin{t}   \\
      0 & - \sin{t} & \cos{t}   
    \end{array}\right),
    \end{equation}
whose inverse is given by 
    
    $$M^{-1}(t) = \left(\begin{array}{ccc}
    p_{y_0}/( p_{y_0}\cos{(t/q)} - y_0 \sin{(t/q)}) &  0  & 0  \\
   0  &  \cos{t} & - \sin{t}   \\  0 & \sin{t} & \cos{t}   
    \end{array}\right).
    $$  
 An easy computation shows that
 $$M^{-1}(0) - M^{-1}(2\pi) = \left(\begin{array}{ccc}
  0 &  0  & 0  \\
   0  &  2\sin^2{(\pi q)} &  \sin{(2\pi q)}   \\
      0 & - \sin{(2\pi q)} & 2\sin^2{(\pi q)}.  
    \end{array}\right).$$
We observe that this matrix has two zeros in the upper right corner of size $1 \times 2$;  the
determinant of the $ 2 \times 2$ matrix which appears in the lower right corner is $4\sin^2{(\pi q)} \neq 0$    
because $q$ is an irrational number. Consequently all the assumptions of Theorem \ref{t3} are satisfied. Therefore 
we must compute the simple zeroes of the function $\CF$ defined in Theorem \ref{t3}.

A straightforward computations shows that 
$$F_1(t,{\bf x}(t,{\bf z})) = \left( -\frac{c (y_0 \cos{(t/q)} + p_{y_0} \sin{(t/q)})^4}
{p_{y_0}\cos{(t/q)} -y_0\sin{(t/q)}}, 0, 0 \right),$$
therefore we have
$$M^{-1}(t) F_1(t,{\bf x}(t,{\bf z}))  =  \left( -\frac{c p_{y_0}(y_0 \cos{(t/q)} + p_{y_0} \sin{(t/q)})^4}{( p_{y_0}\cos{(t/q)} 
-y_0\sin{(t/q)})^2}, 0, 0) \right).$$

Let 
\begin{eqnarray}\label{f1}
f_1(y_0) &=& \frac{1}{2 \pi} \int_0^{2\pi}  \left( -\frac{c p_{y_0}(y_0 \cos{(t/q)} +
 p_{y_0} \sin{(t/q)})^4}{( p_{y_0}\cos{(t/q)} -y_0\sin{(t/q)})^2} \right)dt \nonumber \\
 &=& 3c p_{y_0}( p_{y_0}^2 + y_0^2)/2 = 3ch q^2 \sqrt{2hq - y_0^2}. \nonumber
 \end{eqnarray}
In the last equality we have gotten $p_{x_0}$  from the energy relation $h=\frac{(p_{x_0}^2 + x_0^2)}{2q}$. So the solutions of $f_1=0$ are $y_0= \pm \sqrt{2hq}$,  which are simple zeroes.  On the other hand we can verify that both zeroes generate the same periodic orbit. Doing the rescaling (\ref{rep}) we get the statement {\bf (b)} of Theorem \ref{main1}.

\medskip
Therefore we have proved that for $q$ an irrational number, in every energy level $h>0$, the Hamiltonian system 
(\ref{eq:ham1}) has at least $2$ periodic orbits, so Theorem \ref{main1} holds.
 
 \medskip
\br\label{remarkf}
Using the methods of averaging theory studied in this paper, we could not obtain any periodic orbit for the Hamiltonian system (\ref{eq:ham1}) when $q$ is a rational number. We have tried to get some information in two different ways, using cartesian coordinates as in statement {\bf (a)} and using a modified kind of polar coordinates in two different planes. In the first way we have obtained the variational equations, but unfortunately we could not solve them. In the second way we have obtained that one of the equations 
that we must solve for obtain the periodic solutions is identically zero. \er

\end{document}